# BUBBLES, CONVEXITY AND THE BLACK–SCHOLES EQUATION

By Erik Ekström[1] and Johan Tysk[2]

*Uppsala University*

A bubble is characterized by the presence of an underlying asset whose discounted price process is a strict local martingale under the pricing measure. In such markets, many standard results from option pricing theory do not hold, and in this paper we address some of these issues. In particular, we derive existence and uniqueness results for the Black–Scholes equation, and we provide convexity theory for option pricing and derive related ordering results with respect to volatility. We show that American options are convexity preserving, whereas European options preserve concavity for general payoffs and convexity only for bounded contracts.

**1. Introduction.** Recently, the issue of modeling financial bubbles has attracted the attention of several authors; see, for example, [4, 11, 14, 15] and [21]. We will here adopt the setup of, for instance, [4] and consider markets where the underlying asset, when discounted, follows a strict local martingale under the risk-neutral probability measure. In [14] it is shown that such examples can only exist in a complete market if the Merton no dominance hypothesis is not valid. We will nevertheless use this framework to study the Black–Scholes equation partly because of the interest in its own right, and partly because the problems that arise in this context also arise when one, for instance, considers stochastic volatility models; compare [7, 11] and [13].

In [4] (compare also, for instance, [21]) it is shown that many standard results in option pricing theory fail. For example, put-call parity does not hold, the Black–Scholes equation can have multiple solutions, the price of an American call exceeds that of a European call and the European call

Received March 2008; revised November 2008.

[1]Supported in part by FMB (The Swedish Graduate School in Mathematics and Computing).

[2]Supported in part by the Swedish Research Council.

price is not convex as a function of the stock price. In the present article we provide further insight on option pricing in markets with bubbles. First we study existence and uniqueness theory for the corresponding Black–Scholes equation. It is shown that the option price, given as a risk-neutral expected value, always solves the equation. In general, however, there is an infinite-dimensional space of solutions of linear growth. This is in stark contrast to the standard case when there is a unique solution of polynomial growth. Nevertheless, uniqueness is recovered for contracts of strictly sublinear growth. Next we provide convexity results for markets with bubbles. As mentioned above, and again unlike the case with no bubbles (see, e.g., [8, 12] and [16]), the convexity of a payoff does not necessarily imply the convexity of a European option price. However, we show that American options remain convexity preserving also in the case of bubbles, whereas European options are convexity preserving only for bounded payoffs and concavity preserving for general contracts. Thus, in this respect, prices of American options are more robust then their European counterparts. As a consequence of the convexity/concavity results, we also provide monotonicity results of option prices with respect to volatility.

The paper is organized as follows. In Section 2 we specify the model, and we discuss a motivating example. In Sections 3 and 4 we study existence and uniqueness of solutions to the Black–Scholes equation. Finally, in Sections 5 and 6 we provide convexity theory for the Black–Scholes equation in the presence of bubbles.

**2. The model setup and a motivating example.** We model the stock price process under the risk-neutral probability measure as the solution to the stochastic differential equation

$$dX(t) = \alpha(X(t), t)\, dW, \tag{1}$$

where $\alpha$ is some given function and $W$ is a standard Brownian motion. For the sake of simplicity we let the risk free rate be zero; the case of a deterministic short rate can be treated analogously. We also assume that $x = 0$ is an absorbing state for the price process, that is, if $X(t) = 0$ at some $t$ then $X$ remains at 0 for all times after $t$ as well.

HYPOTHESIS 2.1. *The function $\alpha(x, t)$ is continuous and locally Hölder continuous in space with exponent $1/2$ on $(0, \infty) \times [0, T]$. Moreover, $\alpha(x, t) > 0$ for positive $x$.*

Hypothesis 2.1 ensures the existence of a unique strong solution absorbed at 0. Note that since $X$ is a nonnegative local martingale, it is a supermartingale. Consequently, $X(t)$ does not explode. Given a continuous payoff function $g:[0,\infty) \to [0,\infty)$, the price at time $t$ of a European option that at time $T$ pays the amount $g(X(T))$ is given by

(2) $$u(x,t) = E_{x,t} g(X(T)),$$

where the indices indicate that $X(t) = x$. We thus adhere to classical arbitrage free pricing theory, and let the price be given by this risk neutral expected value even if the corresponding Black–Scholes equation has multiple solutions. Note that if $g$ is of at most linear growth, then $u$ is finite since $X$ is a supermartingale. The corresponding Black–Scholes equation is

(3) $$u_t(x,t) + \tfrac{1}{2}\alpha^2(x,t) u_{xx}(x,t) = 0$$

for $(x,t) \in (0,\infty) \times (0,T)$, with terminal condition

(4) $$u(x,T) = g(x).$$

If the process $X$ hits the boundary $x=0$ with positive probability, then boundary values need to be specified. On the other hand, also if the boundary is hit with probability zero, boundary values are beneficial for instance from a numerical point of view. Recalling that the discount rate is 0 and that the stock price is absorbed at the bankruptcy level $x=0$, we find that the appropriate boundary condition is

(5) $$u(0,t) = g(0).$$

It is well known that if the diffusion coefficient $\alpha$ satisfies a linear growth condition in $x$, that is,

(6) $$|\alpha(x,t)| \leq C(1+x)$$

for some constant $C$, then $u$ defined in (2) is the unique classical solution to (3)–(5) of at most polynomial growth, provided that the payoff function $g$ is of at most polynomial growth (see [17] and the references therein).

In [4] and [11], models in which the underlying asset is a strict local martingale (i.e., a local martingale which is not a martingale) are proposed as models for bubbles. Typical examples include diffusion models in which the linear bound (6) is violated. Another example is furnished by stochastic volatility models, compare [1, 20] and [23]. In this paper we consider models in which (6) does not necessarily hold. In a related article [7] we continue this study to stochastic volatility models.

The following example was studied in [4] and [5], and it shows that the Black–Scholes equation may have multiple solutions. It also shows that convexity of European options is not preserved in general.

EXAMPLE.  Assume that $X$ is given by

$$dX(t) = \sigma X^2(t)\, dW,$$

where $\sigma > 0$ is a constant. Then $X$ is a strict local martingale; see [18] or Theorem 4.1 below. Moreover, the density of $X(T)$ is

$$P(X(T) \in dy) = \frac{x}{y^3\sqrt{2\pi\sigma^2(T-t)}}$$
$$\times \left(\exp\left\{-\frac{(1/y-1/x)^2}{2\sigma^2(T-t)}\right\} - \exp\left\{-\frac{(1/y+1/x)^2}{2\sigma^2(T-t)}\right\}\right),$$

provided $X(t) = x$. It follows that

$$(7) \qquad u(x,t) := E_{x,t}X(T) = x\left(1 - 2\Phi\left(\frac{-1}{x\sigma\sqrt{T-t}}\right)\right)$$

for $t < T$, where $\Phi$ is the distribution function of the standard normal distribution. The corresponding Black–Scholes equation is

$$u_t + \tfrac{1}{2}\sigma^2 x^4 u_{xx} = 0$$

with boundary conditions

$$\begin{cases} u(x,T) = x, \\ u(0,t) = 0. \end{cases}$$

It is straightforward to check that $u$ defined in (7) solves the equation; alternatively, this follows from Theorem 3.2 below. Clearly also $v(x,t) \equiv x$ is a solution to the equation; thus uniqueness of solutions to the Black–Scholes equation fails in general. Moreover,

$$u_{xx}(x,t) = \frac{-2}{x^4\sigma^3(T-t)^{3/2}}\varphi\left(\frac{-1}{x\sigma\sqrt{T-t}}\right),$$

where $\varphi$ is the density of the standard normal distribution. Thus $u$ is strictly concave in $x$ for each fixed $t < T$, so convexity of the payoff function is not preserved. Furthermore, $u$ is decreasing in $\sigma$, which again is the opposite behavior of the standard case.

The above example shows that we are in the somewhat peculiar situation of having a simple and explicit solution to the Black–Scholes equation which is not representing the option price. This is of course a very serious problem if one uses PDE methods to study option prices. A series of natural questions arise:

- Does the option price as given by the stochastic representation always solve the Black–Scholes equation?
- Can one give a general description of models for which uniqueness to the Black–Scholes equation fails, and does uniqueness hold for some restricted class of payoff functions?
- For which contracts is convexity (or concavity) preserved?

- Is preservation of convexity and concavity related to option price monotonicity with respect to volatility?

The first question is answered affirmatively in Section 3, the second one is treated in Section 4 and the remaining questions are investigated in Sections 5 and 6.

**3. Option prices are classical solutions to the Black–Scholes equation.** Since we do not impose any growth restrictions on $\alpha$ at spatial infinity, the Black–Scholes equation (3) falls outside the standard theory for parabolic equations. Indeed, the second-order coefficient $\alpha^2/2$ of the Black–Scholes equation is allowed to grow more than quadratically in the underlying variable, thus violating a standard assumption; compare [9], Theorem 6.4.3. The main result of this section says that the option price $u(x,t)$ given by (2) is indeed a classical solution to (3)–(5); compare [17] for the standard case.

DEFINITION 3.1. A continuous function $v:[0,\infty) \times [0,T] \to \mathbb{R}$ is a *classical solution* to the Black–Scholes equation if $v \in C^{2,1}((0,\infty) \times [0,T))$ and (3)–(5) are satisfied.

THEOREM 3.2. *Assume that the payoff function $g:[0,\infty) \to [0,\infty)$ is continuous and of at most linear growth. Then the option price $u$ defined in (2) is a classical solution to the Black–Scholes equation. Moreover, it is of at most linear growth.*

PROOF. Assume that $g$ satisfies $g(x) \leq C(1+x)$. Then

$$u(x,t) = E_{x,t}g(X(T)) \leq CE_{x,t}(1+X(T)) \leq C(1+x)$$

since $X$ is a supermartingale. Consequently, $u$ is of at most linear growth.

It remains to show that $u$ is a classical solution to the Black–Scholes equation. To indicate the dependence on the starting point, let $X_{x,t}(s)$ be the solution to (1), starting at time $t$ at the point $x$. Moreover, let

$$X_{x,t}^M(s) = X_{x,t}(s \wedge \gamma_M),$$

where

$$\gamma_M = \inf\{s \geq t : X_{x,t}(s) \geq M\}.$$

Then $X^M$ is a bounded local martingale, hence a martingale. Since the paths of $X_{x,t}^M$ and $X_{y,t}^M$, being driven by the same Brownian motion, do not cross each other (see Theorem IX.3.7 in [22]), we have

(8) $$E|X_{x,t}^M(T) - X_{y,t}^M(T)| = |x-y|.$$

For continuity in the initial time variable, let $t_1 \leq t_2$ with $t_2 - t_1 \leq \delta$, and note that

$$(9) \qquad E|X^M_{x,t_1}(t_2) - x|^2 = E\int_{t_1}^{t_2} \alpha^2(X^M_{x,t_1}(s), s)1_{\{s \leq \tau_M\}}\,ds \leq C\delta,$$

where $C = \max\{\alpha^2(x,t) : (x,t) \in [0,M] \times [0,T]\}$. Let $\varepsilon > 0$ and $\eta > 0$. Conditioning on $X^M_{x,t_1}(t_2)$, it follows from (8) that the conditional probability

$$P(|X^M_{x,t_2}(T) - X^M_{x,t_1}(T)| > \varepsilon; |X^M_{x,t_1}(t_2) - x| \leq \delta_0) \leq \eta/2$$

for some $\delta_0 > 0$. Let $\delta = \eta\delta_0^2/(2C)$. Chebyshev's inequality applied to (9) yields

$$P(|X^M_{x,t_1}(t_2) - x| > \delta_0) \leq \frac{C\delta}{\delta_0^2} = \eta/2.$$

Consequently,

$$(10) \qquad P(|X^M_{x,t_2}(T) - X^M_{x,t_1}(T)| > \varepsilon) \leq \eta$$

provided $t_2 - t_1 \leq \delta$. Now (8) and (10) show that $X^M_{y,s}(T)$ tends to $X^M_{x,t}(T)$ in probability as $(y,s) \to (x,t)$.

Next, let $g$ be a nonnegative concave function with $g(0) = 0$. We claim that the function

$$u^M(x,t) = Eg(X_{x,t}(T))1_{\{T < \gamma_M\}}$$

satisfies

$$u^M_t + \frac{\alpha^2}{2}u^M_{xx} = 0 \qquad \text{for } (x,t) \in (0,M) \times [0,T]$$

and

$$\lim_{(y,s) \to (x,T)} u^M(y,s) = g(x) \qquad \text{for all } x \in (0,M).$$

To see this, let $g_n : [0,M] \to \mathbb{R}$ be an increasing sequence of continuous functions such that $g_n(M) = 0$, and such that $g_n(x) \uparrow g(x)$ for each $x \in [0,M)$. Now, each $g_n$ is bounded, so the family

$$\{g_n(X^M_{y,s}(T))\}_{(y,s) \in U},$$

where $U$ is some small neighborhood of $(x,t)$, is uniformly integrable. Consequently, the functions

$$u^M_n(x,t) = Eg_n(X^M_{x,t}(T))$$

are continuous on $[0,M] \times [0,T]$. Since $\alpha(x,t) > 0$ for positive $x$, a continuous stochastic solution is a classical solution, so $u^M_n$ satisfies the Black–Scholes equation in $(0, M) \times [0, T]$; compare, for example, Theorem 2.7 in [17]. By monotone convergence,

$$u_n^M(x, t) \uparrow u^M(x, t) \tag{11}$$

as $n \to \infty$, so it follows from interior Schauder estimates (see [3] or [19]) that also the function $u^M$ solves the Black–Scholes equation on $(0, M) \times [0, T]$. Moreover, (11) yields that

$$\liminf_{(y,s) \to (x,T)} u^M(y, s) \geq g(x).$$

The reverse inequality also holds, since the concavity of $g$ implies that $u^M(x, t) \leq g(x)$ for all $(x, t)$, thus finishing the proof of the claim.

Now, if $u$ is defined as in (2), that is,

$$u(x, t) = E g(X_{x,t}(T)),$$

then $u^M(x, t) \uparrow u(x, t)$ as $M \to \infty$ by monotone convergence. Thus interior Schauder estimates show that $u$ solves the Black–Scholes equation on $(0, \infty) \times [0, T]$. Moreover, the same technique as used above shows that $u$ is continuous up to $t = T$, and since $u \leq g$ we find that $u$ is also continuous up to $x = 0$ with $u(0, t) = g(0) = 0$. This finishes the proof in the case when $g$ is concave and $g(0) = 0$. By linearity, the result also holds for payoffs that are linear combinations of nonnegative concave functions. Note that each smooth function $g$ with

$$\int_0^\infty |g''(y)| \, dy < \infty \tag{12}$$

can be written as a difference of two nonnegative concave functions. Therefore, we approximate $g$ by sequences of smooth functions $g_n$ and $g^n$ that satisfy (12) such that

$$g_n(x) \uparrow g(x)$$

and

$$g^n(x) \downarrow g(x)$$

as $n \to \infty$. Then the result for a general continuous payoff function $g$ follows by an argument involving monotone convergence and Schauder estimates as above. □

**4. Uniqueness and nonuniqueness results for the Black–Scholes equation.** It follows from Theorem 3.2 that uniqueness of solutions of linear growth to the Black–Scholes equation fails if $X$ is a strict local martingale. Indeed, in that case both $u_1(x, t) = E_{x,t} X(T)$ and $u_2(x, t) = x$ solve the Black–Scholes equation with the same boundary conditions; compare the example in Section 2. The following result tells us when the stock price process is a strict local martingale, thus in particular implying that uniqueness to the Black–Scholes equation is lost. For necessary and sufficient conditions for time-homogeneous exponential local martingales to be martingales; see [2].

THEOREM 4.1. *If the volatility coefficient $\alpha$ satisfies*
$$\alpha^2(x,t) \geq \varepsilon x^\eta$$
*for all $(x,t) \in [\varepsilon^{-1}, \infty) \times [0,T]$, where $\varepsilon > 0$ and $\eta > 2$ are constants, then the price process $X$ is a strict local martingale. Moreover, for any time bounded away from expiry, the stock option price is $o(x^\delta)$ for any positive $\delta$, and if $\eta > 3$ then the stock option price is bounded.*

PROOF. Let
$$h(x,t) = e^{M(T-t)} \frac{x}{1 + x^\beta (T-t)^m},$$
where $0 < \beta < 1$. We claim that $m$ and $M$ can be chosen so that $h$ is a supersolution to Black–Scholes equation, that is,

(13) $$h_t + \frac{\alpha^2}{2} h_{xx} \leq 0.$$

To see this, first note that
$$e^{-M(T-t)}(1 + x^\beta(T-t)^m)^3 h_t(x,t)$$
$$= -Mx(1 + x^\beta(T-t)^m)^2 + mx^{\beta+1}(T-t)^{m-1}(1 + x^\beta(T-t)^m)$$

and
$$e^{-M(T-t)}(1 + x^\beta(T-t)^m)^3 h_{xx}(x,t)$$
$$= -\beta(\beta+1)x^{\beta-1}(T-t)^m - \beta(1-\beta)x^{2\beta-1}(T-t)^{2m} \leq 0.$$

Thus, it suffices to choose $m$ and $M$ so that

(14) $$\begin{aligned} & Mx + \tfrac{1}{2}(1+\beta)\beta\alpha^2(x,t)x^{\beta-1}(T-t)^m \\ & \quad + \tfrac{1}{2}\beta(1-\beta)\alpha^2(x,t)x^{2\beta-1}(T-t)^{2m} \\ & \geq mx^{\beta+1}(T-t)^{m-1} + mx^{2\beta+1}(T-t)^{2m-1} \end{aligned}$$

at all points $(x,t) \in [0,\infty) \times [0,T]$.

First we consider large values of $x$, for which by assumption $\alpha^2 \geq \varepsilon x^\eta$. Then it suffices to show

(15) $$\begin{aligned} & Mx + \frac{\varepsilon}{2}(1+\beta)\beta x^{\beta-1+\eta}(T-t)^m + \frac{\varepsilon}{2}\beta(1-\beta)x^{2\beta-1+\eta}(T-t)^{2m} \\ & \geq mx^{\beta+1}(T-t)^{m-1} + mx^{2\beta+1}(T-t)^{2m-1} \end{aligned}$$

at all points $(x,t) \in [\varepsilon^{-1}, \infty) \times [0, T]$. Pick $\eta' \in (2, \eta)$. For $x \geq c_1(m)$, where $c_1$ is a large enough constant depending on $m$, and $T - t \geq x^{2-\eta'}$, the second and third terms on the left-hand side dominate the right-hand side of (15). To see this, we start from the last two terms on left-hand side of (15) and use the lower bound on $T - t$ to obtain

$$\frac{\varepsilon}{2}(1+\beta)\beta x^{\beta-1+\eta}(T-t)^m + \frac{\varepsilon}{2}\beta(1-\beta)x^{2\beta-1+\eta}(T-t)^{2m}$$
$$\geq \frac{\varepsilon}{2}(1+\beta)\beta x^{\beta+1+\eta-\eta'}(T-t)^{m-1} + \frac{\varepsilon}{2}\beta(1-\beta)x^{2\beta+1+\eta-\eta'}(T-t)^{2m-1},$$

which dominates the right-hand side of (15) for large $x$ since $\eta > \eta'$.

Now, for small values of $T - t$, we use the first term on the left-hand side of (15) to dominate the right-hand side, by choosing $m$ sufficiently large and again comparing exponents of $x$. Specifically, for $m > 1 + \beta/(\eta' - 2)$, $x$ dominates the right-hand side for $x \geq c_2(m)$ and $T - t \leq x^{2-\eta'}$, where $c_2$ is chosen sufficiently large. Thus (14) is established for $x \geq c = \max(c_1, c_2)$ and $M \geq 1$.

Next, for the remaining values $x < c$, inequality (14) is established by choosing $M \geq mc^\beta T^{m-1} + mc^{2\beta}T^{2m-1}$. We note that to establish the inequality for small values of $x$, the last two terms on the left-hand side are not needed. Thus we only need to bound the volatility coefficients from below for large values of $x$ as is done in the statement of the theorem.

Now, since $h$ is nonnegative, it follows from Itô's formula and (13) that the process $h(X(s), s)$ is a supermartingale. Hence, for $t < T$,

$$x > h(x, t) \geq E_{x,t}h(X(T), T) = E_{x,t}X(T),$$

which implies that $X$ is a strict local martingale.

Finally, if $\eta > 3$, then $\beta$ can be chosen to be 1, and $m$ and $M$ can be chosen so that the two first terms on the left-hand side of (15) together dominate the terms on the right. $\square$

REMARK. We note that it is not enough to assume that $\frac{\alpha^2(x,t)}{x^2} \to \infty$ for $X$ to be a strict local martingale. If $\alpha^2(x) = x^2 \ln x$, then the price of the stock option is $x$, implying that $X$ is a martingale. This can be seen by considering the supersolution $h(x,t) = e^{M(T-t)}(1 + x \ln x)$ and arguing as in the proof of Theorem 4.3 below.

REMARK. If $X$ is a strict local martingale, then the space of classical solutions to the Black–Scholes equation is infinite dimensional. Indeed, let $\tilde{T} \leq T$ and define

$$(16) \qquad v_{\tilde{T}}(x,t) = \begin{cases} x - E_{t,x}X(\tilde{T}), & \text{for } t \leq \tilde{T} \leq T, \\ 0, & \text{for } \tilde{T} \leq t \leq T. \end{cases}$$

Using boundary Schauder estimates (see [10]), it follows that $v_{\tilde{T}}$ is a classical solution to the homogeneous equation on $[0, \infty) \times [0, T]$ (technically speaking, an additional Hölder continuity of $\alpha$ in the time variable needs to be imposed to apply these estimates). The set $\{\lambda v_{\tilde{T}} : \lambda \in \mathbb{R} \text{ and } \tilde{T} \in (0, T]\}$ is therefore an infinite-dimensional space of solutions to the homogeneous equation.

To get uniqueness of solutions, one needs to narrow the class of considered functions. It turns out that the appropriate class in which uniqueness holds is the class of functions of strictly sublinear growth.

DEFINITION 4.2. A function $f : [0, \infty) \to \mathbb{R}$ is of strictly sublinear growth if $\lim_{x \to \infty} \frac{f(x)}{x} = 0$. A function $v : [0, \infty) \times [0, T] \to \mathbb{R}$ is of strictly sublinear growth if $|v(x,t)| \leq f(x)$ for some $f : [0, \infty) \to [0, \infty)$ of strictly sublinear growth.

THEOREM 4.3. *Assume that the payoff function $g$ is of strictly sublinear growth. Then the option price $u$ is the unique classical solution of strictly sublinear growth to the Black–Scholes equation.*

PROOF. It follows from Theorem 3.2 that $u$ is a classical solution to the Black–Scholes equation. Moreover, if $g$ is of strictly sublinear growth, then so is its smallest concave majorant $\overline{g}$. Consequently,

$$u(x,t) \leq E_{x,t}\overline{g}(X(T)) \leq \overline{g}(E_{x,t}X(T)) \leq \overline{g}(x),$$

where we use Jensen's inequality and the monotonicity of $\overline{g}$. Thus $u$ is a solution of strictly sublinear growth.

To prove uniqueness, assume that $v$ is a classical solution to

$$\begin{cases} v_t = \frac{1}{2}\alpha^2 v_{xx}, \\ v(x,0) = 0, \\ v(0,t) = 0 \end{cases}$$

of strictly sublinear growth (for simplicity, we have performed a standard change of variables $t \to T - t$, and the terminal condition is then replaced by an initial condition). Define

$$h(x,t) = e^t(1 + x)$$

and let $v^\varepsilon(x,t) = v(x,t) + \varepsilon h(x,t)$. Then

(17) $$v_t^\varepsilon - \tfrac{1}{2}\alpha^2 v_{xx}^\varepsilon = \varepsilon h_t - \varepsilon \tfrac{1}{2}\alpha^2 h_{xx} = \varepsilon e^t(1+x) > 0.$$

Next, define

$$\Gamma = \{(x,t) \in [0, \infty) \times [0, T] : v^\varepsilon(x,t) < 0\}.$$

Assuming that $\Gamma \neq \varnothing$, let

(18) $\qquad t_0 = \inf\{t \in [0,T] : (x,t) \in \overline{\Gamma} \text{ for some } x \in [0,\infty)\}.$

Since $v^\varepsilon$ is of linear growth, the set $\Gamma$ is contained in $[0,M] \times [0,T]$ for some constant $M$. Consequently, $\overline{\Gamma}$ is compact, and the infimum in (18) is attained for some point $(x_0, t_0)$ with $v^\varepsilon(x_0, t_0) = 0$. Since $v^\varepsilon(0,t) = \varepsilon e^t > 0$ and $v^\varepsilon(x,0) = 1 + x > 0$, we have $x_0 > 0$ and $t_0 > 0$. Now, by the definition of $t_0$, the function $x \mapsto v^\varepsilon(x, t_0)$ has a minimum at $x = x_0$. Thus $v^\varepsilon_{xx} \geq 0$. Similarly, $v^\varepsilon_t \leq 0$, so

$$v^\varepsilon_t - \tfrac{1}{2}\alpha^2 v^\varepsilon_{xx} \leq 0,$$

which contradicts (17). This contradiction shows that $\Gamma$ is empty, so $v^\varepsilon \geq 0$. Since $\varepsilon$ is arbitrary, $v \geq 0$ follows. Finally, the same argument applied to $-v$ shows the reverse inequality, so $v \equiv 0$, which proves uniqueness. $\square$

In the case of multiple solutions of the Black–Scholes equation, the solution $u$ in (2) can be characterized as the smallest nonnegative supersolution; see [11]. Indeed, it follows from the Itô formula that the process $v(X(t),t)$ is a supermartingale, provided $v$ is a nonnegative supersolution. Consequently,

$$v(x,t) \geq E_{x,t} v(X(T), T) \geq E_{x,t} g(X(T)) = u(x,t).$$

In addition to this characterization, we can identify $u$ as the limit of a sequence of solutions to the Black–Scholes equation with bounded payoffs.

PROPOSITION 4.4. *Let the payoff function $g$ be of at most linear growth. Then the option price $u$ in (2) is the limit of the (unique) classical solution to the Black–Scholes equation with terminal value $\min(g, M)$ as $M$ tends to infinity.*

PROOF. By dominated convergence,

$$u(x,t) = E_{x,t}\Big[\lim_{M \to \infty} g(X(T)) \wedge M\Big] = \lim_{M \to \infty} E_{x,t}[g(X(T)) \wedge M]. \quad \square$$

**5. Convexity theory for European options.** In standard Black–Scholes theory, it is well known that prices of options with convex payoffs are convex in the stock price; compare [8, 12] and [16]. However, this result is not true for markets with bubbles. Indeed, as noted in Section 2, the function $u$ defined in (7) is strictly concave in $x$.

In this section we show that models for bubbles are convexity preserving for *bounded* contracts, and concavity preserving for *all* contracts. The lack of symmetry is due to the fact that we consider only nonnegative payoffs. As in the standard case, preservation of convexity or concavity implies monotonicity properties with respect to volatility. To formulate this, assume that $\alpha_1$ and $\alpha_2$ are two nonnegative volatility functions satisfying

$$\alpha_1(x,t) \leq \alpha_2(x,t)$$

for all $(x,t) \in [0,\infty) \times [0,T]$, and let $u_1(x,t)$ and $u_2(x,t)$ be the corresponding option prices.

THEOREM 5.1. *Assume that $g$ is concave. Then $u(x,t)$ is concave in $x$ for any $t \in [0,T]$. Moreover, the option price is decreasing in the volatility, that is, $u_1(x,t) \geq u_2(x,t)$ for all $(x,t) \in [0,\infty) \times [0,T]$.*

*Similarly, if $g$ is convex and bounded, then $u(x,t)$ is convex in $x$ for any $t \in [0,T]$. Moreover, the option price is increasing in the volatility, that is, $u_1(x,t) \leq u_2(x,t)$ for all $(x,t) \in [0,\infty) \times [0,T]$.*

REMARK. In fact, the result holds for concave payoff functions that are bounded from below, and convex functions bounded from above, respectively. The asymmetry between the conditions on the payoff functions in the statement of the theorem is due to our assumption that $g$ is nonnegative.

PROOF OF THEOREM 5.1. Assume first that $g$ is bounded and concave. Without loss of generality we also assume that $g(0) = 0$. Let $u_K$ be the corresponding option price where the volatility is given by $\alpha(x,t) \wedge K$. It is well known that $u_K$ is concave; see, for instance, [16]. It also follows from this reference that this sequence of functions is decreasing in $K$. We let

$$\overline{u}(x,t) = \lim_{K \to \infty} u_K(x,t)$$

denote its limit. By interior Schauder estimates (see [3] and [19]), $\overline{u}$ solves the Black–Scholes equation at all points in $(0,\infty) \times [0,T)$. Also $\overline{u}(0,t) = 0$ since this holds for each function in the decreasing sequence $u_K$ of nonnegative functions. Now, for any positive number $b$, consider a continuous function $\psi \leq 1$ that is identically 1 in a neighborhood of $b$ and with support in the interval $[\frac{b}{2}, 2b]$. Let $v_{K,b}$ be the solution of the pricing equation with volatility $\alpha(x,t) \wedge K$ and contract function $g\psi$, but with vanishing Dirichlet data at $x = \frac{b}{2}$ and $x = 2b$. We note that $v_{K,b}$ does not depend on $K$ if $K$ is large enough. Let $v_b(x,t) = \lim_{K \to \infty} v_{K,b}(x,t)$. By the maximum principle, $v_{K,b} \leq u_K$ and thus $v_b \leq \overline{u}$. Hence $v_b \leq \overline{u} \leq u_K$ for any $K$. But $v_b$ is known to be continuous up to the boundary $t = T$ by classical theory and $u_K$ by [17]. Hence $\overline{u}$ is continuous up to the boundary with boundary values given by $g$ in a neighborhood of $b$. Carrying out the same approximation argument we conclude that this holds at any point. Hence $\overline{u}$ is a classical solution to the Black–Scholes equation with $g$ as terminal value. By the uniqueness result Theorem 4.3 it follows that $\overline{u} = u$. Since $\overline{u}$ is the limit of a sequence of concave functions, it is concave itself, so the concavity of $u$ follows.

Next, if $g$ is concave but not bounded, then let $g^M = g \wedge M$, and let $u^M$ be the corresponding option price. It follows from above that $u^M$ is concave in $x$. Consequently, it follows from Proposition 4.4 that $u$ is concave.

Finally, the monotonicity in volatility stated in the theorem follows from the approximation argument above and the corresponding monotonicity in the bubble-free case; see for instance [16]. The second part of the theorem follows by a similar argument. □

REMARK. The crucial ingredient in the above proof is the uniqueness result Theorem 4.3. In the absence of this result, for instance for general convex functions of linear growth, the limiting function $\overline{u}$ will still be a solution to the Black–Scholes equation, but it will no longer be the solution given by the stochastic representation.

**6. Convexity theory for American options.** We have seen above that European options in general do not preserve convexity. In this section we show that American option prices are more well behaved in this respect.

When pricing American options one cannot without loss of generality assume that the short rate is 0. Therefore, let $X$ be the solution to

$$dX = rX(t)\,dt + \alpha(X(t), t)\,dW,$$

and define the American option price $U : [0, \infty) \times [0, T]$ by

$$U(x, t) = \sup_{t \leq \tau \leq T} E_{x,t} e^{-r(\tau - t)} g(X(\tau)).$$

Here the supremum is over random times $\tau$ that are stopping times with respect to the completion of the filtration generated by $W$.

THEOREM 6.1. *Assume that the payoff function $g$ is convex and of at most linear growth. Then the American option price $U(x, t)$ is convex in $x$. Moreover, $U$ is increasing in the volatility.*

PROOF. For notational simplicity we let $t = 0$. Since $g$ is convex, nonnegative and of at most linear growth, it is either decreasing to a limit $\lim_{x \to \infty} g(x) = constant \geq 0$ or it satisfies $\lim_{x \to \infty} (g(x)/x) = constant > 0$.

First we treat the case $\lim_{x \to \infty} (g(x)/x) = \gamma > 0$. For simplicity, assume that $\gamma = 1$, so that $g$ satisfies

$$g(x) \leq C + x \tag{19}$$

for some $C > 0$. For $M > 0$, let $X^M$ be the solution to

$$dX^M(t) = rX^M(t)\,dt + \alpha^M(X^M(t), t)\,dW,$$

where $\alpha^M(x,t) = \alpha(x,t) \wedge M$. For $\varepsilon > 0$, define

$$U_\varepsilon^M(x,t) = \sup_{t \leq \tau \leq T} E_{x,t} e^{-r(\tau-t)} g_\varepsilon(X^M(\tau)),$$

where $g_\varepsilon(x) = (1-\varepsilon)g(x)$. It follows from (19) that

$$U_\varepsilon^M(x,0) \leq (1-\varepsilon)(C+x).$$

Consequently, for each fixed $\varepsilon > 0$, there exists $x_0$ such that

$$U_\varepsilon^M(x,t) \leq g(x) \leq U(x,t)$$

for all $x \geq x_0$ and all $t \in [0,T]$. If $M$ is large enough, then $X^M(t)$ and $X(t)$ coincide for $t \leq \tau_0 = \inf\{s : X(s) \geq x_0\}$, so it follows from the strong Markov property that $U_\varepsilon^M(x,t) \leq U(x,t)$ for $x \leq x_0$. Consequently,

(20) $$\lim_{\varepsilon \to 0} \lim_{M \to \infty} U_\varepsilon^M \leq U$$

(both limits exist since $U_\varepsilon^M$ is increasing in $M$—see [6, 8] and [12]—and decreasing in $\varepsilon$). To demonstrate the reverse inequality, let $\tau^*$ be a $\delta$-optimal stopping time for $U$, that is,

$$U(x,0) \leq E_{x,0} e^{-r\tau^*} g(X(\tau^*)) + \delta.$$

Then

$$\lim_{M \to \infty} U_\varepsilon^M(x,0) \geq \liminf_{M \to \infty} E_{x,0} e^{-r\tau^*} g_\varepsilon(X^M(\tau^*))$$

$$\geq (1-\varepsilon) E_{x,0} e^{-r\tau^*} \liminf_{M \to \infty} g(X^M(\tau^*))$$

$$= (1-\varepsilon) E_{x,0} e^{-r\tau^*} g(X(\tau^*))$$

$$\geq (1-\varepsilon)(U(x,0) - \delta),$$

where we used Fatou's lemma and the fact that $X^M(\tau^*) \to X(\tau^*)$ almost surely as $M \to \infty$. Since $\delta$ is arbitrary, we find that

$$\lim_{\varepsilon \to 0} \lim_{M \to \infty} U_\varepsilon^M \geq U.$$

Now, recall that $U_\varepsilon^M$ is convex in $x$; see [6] or [8]. Since the pointwise limit of a sequence of convex functions is convex, it follows that $U$ is convex in $x$. Similarly, $U_\varepsilon^M$ is increasing in the volatility (see [6] or [8]), so $U$ is also increasing in the volatility.

Next, consider the case of a decreasing payoff function satisfying

$$\lim_{x \to \infty} g(x) = \gamma \geq 0.$$

In this case, define

$$U^M(x,0) = \sup_{\tau \leq \tau_M \wedge T} E_{x,0} e^{-r\tau} g(X(\tau)),$$

where
$$\tau_M = \inf\{t : X(t) \geq M\}$$
is the first passage time of $X$ over the level $M$. Clearly, $U^M \leq U$ and $U^M$ is increasing in $M$, so
$$\lim_{M \to \infty} U^M(x,0) \leq U(x,0).$$
The reverse inequality follows by the following argument. Take $\delta > 0$, and let $\tau^*$ be a $\delta$-optimal stopping time for $U$, that is,
$$U(x,0) \leq E_{x,0} e^{-r\tau^*} g(X(\tau^*)) + \delta.$$
Then
$$\begin{aligned}
\lim_{M \to \infty} U^M(x,0) &\geq \liminf_{M \to \infty} E_{x,0} e^{-r(\tau^* \wedge \tau_M)} g(X(\tau^* \wedge \tau_M)) \\
&\geq \liminf_{M \to \infty} E_{x,0} e^{-r\tau^*} g(X(\tau^*)) 1_{\{\tau^* \leq \tau_M\}} \\
&\geq E_{x,0} e^{-r\tau^*} g(X(\tau^*)) \\
&\geq U(x,0) - \delta,
\end{aligned}$$
where the third inequality follows from Fatou's lemma and the fact that $\gamma_M \to \infty$ as $M \to \infty$. Since $\delta > 0$ is arbitrary, this finishes the proof of $\lim_{M \to \infty} U^M = U$.

It is straightforward to check that $U^M(x,t)$ is convex in $x$ and increasing in the volatility. Indeed, this can be shown by approximating $U^M$ with Bermudan options, all of which are convex in $x$ and increasing in the volatility (compare [6]). Both these properties are therefore inherited by $U$, which finishes the proof. □

REMARK. The above proof in the case of a decreasing convex payoff function does not carry over to the general case since
$$U^M(x,0) = \sup_{\tau \leq \tau_M} E_{x,0} e^{-r\tau} g(X(\tau))$$
is not convex in general.

REMARK. One may note that American options are trivially concavity preserving. Indeed, if $g$ is concave, then $e^{-rt} g(X(t))$ is a supermartingale, thus implying that $V(x,t) = g(x)$ (see [6] for the case when $\alpha$ satisfies (6)).

**Acknowledgments.** We thank the anonymous referees for a thorough reading and valuable comments.

## REFERENCES

[1] ANDERSEN, L. B. G. and PITERBARG, V. V. (2007). Moment explosions in stochastic volatility models. *Finance Stoch.* **11** 29–50. MR2284011

[2] BLEI, S. and ENGELBERT, H.-J. (2008). On exponential local martingales associated with strong Markov continuous local martingales. *Stoch. Proc. Appl.* To appear.

[3] BRANDT, A. (1969). Interior Schauder estimates for parabolic differential- (or difference-) equations via the maximum principle. *Israel J. Math.* **7** 254–262. MR0249803

[4] COX, A. M. G. and HOBSON, D. G. (2005). Local martingales, bubbles and option prices. *Finance Stoch.* **9** 477–492. MR2213778

[5] DELBAEN, F. and SCHACHERMAYER, W. (1994). Arbitrage and free lunch with bounded risk for unbounded continuous processes. *Math. Finance* **4** 343–348.

[6] EKSTRÖM, E. (2004). Properties of American option prices. *Stochastic Process. Appl.* **114** 265–278. MR2101244

[7] EKSTRÖM, E. and TYSK, J. (2008). The Black–Scholes equation in stochastic volatility models. Preprint.

[8] EL KAROUI, N., JEANBLANC-PICQUÉ, M. and SHREVE, S. E. (1998). Robustness of the Black and Scholes formula. *Math. Finance* **8** 93–126. MR1609962

[9] FRIEDMAN, A. (1976). *Stochastic Differential Equations and Applications. Vols. 1 and 2. Probability and Mathematical Statistics* **28**. Academic Press, New York. MR0494490

[10] FRIEDMAN, A. (1958). Boundary estimates for second order parabolic equations and their applications. *J. Math. Mech.* **7** 771–791. MR0108648

[11] HESTON, S., LOEWENSTEIN, M. and WILLARD, G. (2007). Options and bubbles. *Rev. Financ. Stud.* **20** 359–390.

[12] HOBSON, D. G. (1998). Volatility misspecification, option pricing and superreplication via coupling. *Ann. Appl. Probab.* **8** 193–205. MR1620358

[13] HOBSON, D. (2008). Comparison results for stochastic volatility models via coupling. *Finance Stoch.* To appear.

[14] JARROW, R. A., PROTTER, P. and SHIMBO, K. (2006). Asset price bubbles in complete markets. In *Advances in Mathematical Finance* 105–130. Birkhäuser, Boston. MR2359365

[15] JARROW, R. A., PROTTER, P. and SHIMBO, K. (2007). Asset price bubbles in incomplete markets. Preprint.

[16] JANSON, S. and TYSK, J. (2003). Volatility time and properties of option prices. *Ann. Appl. Probab.* **13** 890–913. MR1994040

[17] JANSON, S. and TYSK, J. (2006). Feynman–Kac formulas for Black–Scholes-type operators. *Bull. London Math. Soc.* **38** 269–282. MR2214479

[18] JOHNSON, G. and HELMS, L. L. (1963). Class $D$ supermartingales. *Bull. Amer. Math. Soc.* **69** 59–62. MR0142148

[19] KNERR, B. F. (1980/81). Parabolic interior Schauder estimates by the maximum principle. *Arch. Rational Mech. Anal.* **75** 51–58. MR592103

[20] LEWIS, A. L. (2000). *Option Valuation Under Stochastic Volatility: With Mathematica Code.* Finance Press, Newport Beach, CA. MR1742310

[21] PAL, S. and PROTTER, P. (2008). Strict local martingales, bubbles and no early exercise. Preprint.

[22] REVUZ, D. and YOR, M. (1999). *Continuous Martingales and Brownian Motion*, 3rd ed. *Grundlehren der Mathematischen Wissenschaften [Fundamental Principles of Mathematical Sciences]* **293**. Springer, Berlin. MR1725357

[23] Sin, C. A. (1998). Complications with stochastic volatility models. *Adv. in Appl. Probab.* **30** 256–268. MR1618849

Department of Mathematics
Uppsala University
Box 480
SE-751 06 Uppsala
Sweden
E-mail: Erik.Ekstrom@math.uu.se
Johan.Tysk@math.uu.se